\documentclass[11pt,a4paper]{article}
\usepackage{a4wide}
\usepackage{amsmath,amssymb}
\usepackage[english]{babel}
\usepackage{graphicx}
\usepackage[numbers]{natbib}
\usepackage{theorem}
\usepackage{titlesec}
\usepackage{ifpdf}
\usepackage{hyperref}
\usepackage{setspace}

\hypersetup{
    bookmarksopen=false,
    bookmarksnumbered=true,
    pdftitle={Mixed Gated/Exhaustive Service in a Polling Model with Priorities},
    pdfauthor={M.A.A. Boon, I.J.B.F. Adan}
}
\ifpdf
  \hypersetup{colorlinks=true,linkcolor=black,urlcolor=black,citecolor=black,pdfpagemode=UseOutlines,plainpages=false,pdfpagelabels}
\else
  \hypersetup{colorlinks=false}
\fi

\titleformat{\subsubsection}[hang]{\normalfont\em}{}{1ex}{}

\theorembodyfont{\upshape}
\theoremheaderfont{\bfseries}

\newtheorem{theorem}{Theorem}[section]
\newtheorem{property}[theorem]{Property}

\numberwithin{equation}{section}

\newcommand{\ee}{\textrm{e}}

\newcommand{\z}{\mathbf{z}}
\providecommand{\href}[2]{#2}

\title{Mixed Gated/Exhaustive Service in a Polling Model with Priorities\footnote{The research was done in the framework of the BSIK/BRICKS project, and of the European Network of Excellence Euro-FGI.}}
\author{M.A.A. Boon\footnote{\textsc{Eurandom} and Department of Mathematics and Computer Science, Eindhoven University of Technology, P.O. Box 513, 5600MB Eindhoven, The Netherlands}\\\href{mailto:marko@win.tue.nl}{marko@win.tue.nl} \and I.J.B.F. Adan\footnotemark[2]\\\href{mailto:iadan@win.tue.nl}{iadan@win.tue.nl} }

\date{February, 2009}

\begin{document}
\maketitle

\begin{abstract}
In this paper we consider a single-server polling system with switch-over times.
We introduce a new service discipline, mixed gated/exhaustive service, that can be used for queues with two types of customers: high and low priority customers. At the beginning of a visit of the server to such a queue, a gate is set behind all customers. High priority customers receive priority in the sense that they are always served before any low priority customers. But high priority customers have a second advantage over low priority customers. Low priority customers are served according to the gated service discipline, i.e. only customers standing in front of the gate are served during this visit. In contrast, high priority customers arriving during the visit period of the queue are allowed to pass the gate and all low priority customers before the gate.

We study the cycle time distribution, the waiting time distributions for each customer type, the joint queue length distribution of all priority classes at all queues at polling epochs, and the steady-state marginal queue length distributions for each customer type. Through numerical examples we illustrate that the mixed gated/exhaustive service discipline can significantly decrease waiting times of high priority jobs. In many cases there is a minimal negative impact on the waiting times of low priority customers but, remarkably, it turns out that in polling systems with larger switch-over times there can be even a positive impact on the waiting times of low priority customers.

\end{abstract}

\section{Introduction}\label{applicationssection}

There are three ways in which one can introduce prioritisation into a polling model. The first type of priority is by changing the server routing such that certain queues are visited more frequently than other queues \cite{boxmaweststrate89,srinivasan91}. This type of prioritisation is quite common in wireless network protocols. A second type of prioritisation is through differentiation of the number of customers that are served during each visit to a queue. This type of prioritisation is inflicted through the usage of different service disciplines. For example, one can serve all customers in a queue before switching to the next queue (exhaustive service), or one can limit the amount of customers that are served to, e.g., only those customers present at the arrival of the server at the queue (gated service). Typically, this will have a negative impact on the waiting times of the customers in queues that are not served exhaustively. The third way of introducing priorities is by changing the order in which customers are served within a queue, which is a popular technique to improve performance of production systems, cf. \cite{boonadanboxma2008,wierman07}. The present paper introduces a new service discipline, referred to as mixed gated/exhaustive service, that combines the last two types of prioritisation.

In the polling model considered in the present paper a single server visits $N$ queues in a fixed, cyclic order. Some, or even all, of the queues contain two types of customers: high and low priority customers. For these queues we introduce a new service discipline, called mixed gated/exhaustive service based on the priority level of the customer. A polling system with high and low priority customers in a queue with purely gated or exhaustive service has been studied in \cite{boonadanboxma2queues2008,boonadanboxma2008}. The mixed gated/exhaustive service discipline can be considered as a mixture of these two service disciplines where low priority customers receive gated service and high priority customers receive exhaustive service. A more detailed description is given in Section \ref{modeldescription}. Since the number of customers served during one visit in a queue with gated service is different from the number served during a visit with exhaustive service, the mixed gated/exhaustive service discipline introduced in the present paper combines the second and the third type of prioritisation. A variation of the model under consideration, namely a polling system where low priority customers are served only if there are no high priority customers present \emph{in any of the queues}, has been studied in \cite{gianinimanfield88}.

Polling models have been studied for many years and because of their practical relevance many papers on polling systems have been written in a mixture of application areas. The survey of Takagi \cite{takagi91} on polling systems and their applications from 1991 is still very valuable, although the last couple of years interest in polling models has revived, partly triggered by many new applications. The motivation for the present paper is to present a service discipline that combines the benefits of the gated and exhaustive service disciplines for priority polling models. The specific application that attracted our attention is in the field of logistics. Consider a make-to-order production system with a single production capacity for multiple products. In many firms encountering this situation, the products are produced according to a fixed production sequence. The production capacity, where the production orders queue up, can be represented as a polling model by identifying each product with a queue and the demand process of a product with the arrival process at the corresponding queue. For a more detailed description of fixed-sequence strategies in the context of make-to-stock production situations, see  \cite{winandsPhD}. In the context of this production setting, the situation with two or more priority levels - as studied in detail in the present paper - is oftentimes encountered in practice, where production departments have to supply both internal and external customers, the latter of which is commonly given a preferential treatment. A different application stems from production scheduling in flexible manufacturing systems where part types are often grouped with other types sharing (almost) similar characteristics, such that no change of machine configuration, i.e. setup time, is required when switching between these part types (see, e.g., \cite{sharafali04}). Since no setup time is required to switch between these  types, it can be seen as customers of different types being served in the same queue. The introduction of priorities can be useful to efficiently differentiate between different parts grouped within one queue. These two applications make the practical relevance of the inclusion of multiple priority levels in the studied polling model evident. Finally, we should keep in the back of our mind that the results of the present paper are certainly not limited to these production settings, but may be used in many other fields where polling models arise, such as communication, transportation and health care (e.g., surgery procedures where an urgency parameter is assigned to each patient).

The present paper is structured as follows: first we discuss the model in more detail and we determine the generating functions (GFs) of the joint queue length distribution of all customers at visit beginnings and completions of each queue. In Section \ref{cycletimesection} we determine the Laplace-Stieltjes Transforms (LSTs) of the distributions of the cycle time, visit times and intervisit times. These distributions are used to determine the marginal queue length distributions and waiting time distributions of high and low priority customers in all queues. The LST of the waiting time distribution is used to compute the mean waiting time of each customer type. A pseudo-conservation law for these mean waiting times is presented in Section \ref{pseudoconservationlawprioritiessection}. Furthermore, we introduce some numerical examples to illustrate typical features of a polling model with mixed gated/exhaustive service. Finally we discuss possible extensions and future research on the topic.

\section{Notation and model description}\label{modeldescription}

The model considered in the present paper is a polling model which consists of $N$ queues, labelled $Q_1, \dots, Q_N$. Throughout the whole paper all indices are modulo $N$, so $Q_{N+1}$ stands for $Q_1$. The queues are visited by one server in a fixed, cyclic order: $1, 2, \dots, N, 1, 2, \dots$. The switch-over time of the server from $Q_i$ to $Q_{i+1}$ is denoted by $S_i$ with LST $\sigma_i(\cdot)$. We assume that all switch-over times are independent and at least one switch-over time is strictly greater than zero. Each queue contains two customer types: high and low priority customers, although the analysis allows any number (greater than zero) of customer types per queue. High priority customers in $Q_i$ are called type $iH$ customers and low priority customers in $Q_i$ are called $iL$ customers, $i=1,\dots,N$. Type $iH$ customers arrive at $Q_i$ according to a Poisson process with intensity $\lambda_{iH}$, and type $iL$ customers arrive at $Q_i$ according to a Poisson process with intensity $\lambda_{iL}$. The service times of type $iH$ and $iL$ customers are denoted by $B_{iH}$ and $B_{iL}$, with LSTs $\beta_{iH}(\cdot)$ and $\beta_{iL}(\cdot)$. All service times are assumed to be independent. We introduce the notation $\rho_{iH} = \lambda_{iH}E(B_{iH})$ and similarly  $\rho_{iL} = \lambda_{iL}E(B_{iL})$. The total occupation rate of the system is $\rho = \sum_{i=1}^N\rho_{i}$, where $\rho_i = \rho_{iH}+\rho_{iL}$ is the fraction of time that the server visits $Q_i$. Service of the customers is gated for low priority customers and exhaustive for high priority customers. In more detail: each queue actually contains two lines of waiting customers: one for the low priority customers and one for the high priority customers. At the beginning of a visit to $Q_i$, a gate is set behind the low priority customers to mark them eligible for service. High priority customers are always served exhaustively until no high priority customer is present. When no high priority customers are present in the queue, the low priority customers standing in front of the gate are served in order of arrival, but whenever a high priority customer enters the queue, he is served before any waiting low priority customers. Service is non-preemptive though, implying that service of a type $iL$ customer is not interrupted by an arriving type $iH$ customer. The visit to $Q_i$ ends when all type $iL$ customers present at the beginning of this visit are served and no high priority customers are present in the queue. Notice that if the arrival intensity $\lambda_{iH}$ equals $0$, then $Q_i$ is served completely according to the gated service discipline. Similarly we can set $\lambda_{iL} = 0$ to obtain a purely exhaustively served queue. Both the gated and the exhaustive service discipline fall into the category of branching-type service disciplines. These are service disciplines that satisfy the following property, introduced by Resing \cite{resing93} and Fuhrmann \cite{fuhrmann81}.

\begin{property}\label{resingproperty}
If the server arrives at $Q_i$ to find $k_i$ customers there, then during the course of the server's visit, each of these $k_i$ customers will effectively be replaced in an i.i.d. manner by a random population having probability generating function $h_i(z_1,\dots,z_N)$, which can be any $N$-dimensional probability generating function.
\end{property}

If $Q_i$ receives gated service, we have $h_i(z_1,\dots,z_N) = \beta_i\left(\sum_{j=1}^N \lambda_j(1-z_j)\right)$, where $\beta_i(\cdot)$ denotes the service time LST of an arbitrary customer in $Q_i$, and $\lambda_i$ denotes his arrival rate. For exhaustive service $h_i(z_1,\dots,z_N) = \pi_i\left(\sum_{j\neq i} \lambda_j(1-z_j)\right)$, where $\pi_i(\cdot)$ is the LST of a busy period distribution in an $M/G/1$ system with only type $i$ customers, so it is the root in $(0,1]$ of the equation $\pi_i(\omega) = \beta_i(\omega + \lambda_i(1 - \pi_i(\omega)))$, $\omega \geq 0$ (cf. \cite{cohen82}, p. 250).

Property \ref{resingproperty} is not satisfied if $Q_i$ receives mixed gated/exhaustive service, because the random population that replaces each of these customers depends on the priority level. In the next section
we circumvent this problem by splitting each queue into two virtual queues, each of which has a branching-type service discipline. This equivalent polling system satisfies Property \ref{resingproperty}, so we can still use the methodology described in \cite{resing93} to find, e.g., the joint queue length distribution at visit beginnings and completions. All other probability distributions that are derived in the present paper can be expressed in terms of (one of) these joint queue length distributions.

\section{Joint queue length distribution at polling epochs}\label{joint}

In the present section we analyse a polling system with all queues having two priority levels and receiving mixed gated/exhaustive service, but in fact each queue would be allowed to have any branching-type service discipline.
Denote the GF of the joint queue length distribution of type $1H, 1L, \dots, NH, NL$ customers at the beginning and the completion of a visit to $Q_i$ by respectively $V_{b_i}(z_{1H}, z_{1L}, \dots, z_{NH}, z_{NL})$ and $V_{c_i}(z_{1H}, z_{1L}, \dots, z_{NH}, z_{NL})$. As discussed in the previous section, the polling model under consideration does not satisfy Property \ref{resingproperty}, which often means that an exact analysis is difficult or even impossible. For this reason we introduce a different polling system that does satisfy Property \ref{resingproperty} and has the same joint queue length distribution at visit beginnings and endings. The equivalent system contains $2N$ queues, denoted by $Q_{1H^*}, Q_{1L^*}, \dots, Q_{NH^*}, Q_{NL^*}$. The switch-over times $S_i$, $i=1,\dots,N$, are incurred when the server switches from $Q_{iL^*}$ to $Q_{(i+1)H^*}$; there are no switch-over times between $Q_{iH^*}$ and $Q_{iL^*}$. Customers in this system are so-called ``smart customers'', introduced in \cite{smartcustomers}, meaning that the arrival rate of each customer type depends on the location of the server. Type $iH^*$ customers arrive in $Q_{iH^*}$ according to arrival rate $\lambda_{iH}$ \emph{unless the server is serving $Q_{iL^*}$}. When the server is serving $Q_{iL^*}$, the arrival rate of type $iH^*$ customers is 0. The reason for this is that we incorporate the service times of all type $iH$ customers that would have arrived during the service of a type $iL$ customer, in the original polling model, into the service time of a type $iL^*$ customer. In our alternative system, type $iL^*$ customers arrive with intensity $\lambda_{iL}$ and have service requirement $B_{iL}^*$ with LST $\beta_{iL}^*(\cdot)$. There is a simple relation between $B_{iL}$ and $B_{iL}^*$, expressed in terms of the LST:
\begin{equation}
\beta_{iL}^*(\omega) = \beta_{iL}(\omega+\lambda_{iH}(1-\pi_{iH}(\omega))).\label{extendedservicetimeL}
\end{equation}
$B_{iL}^*$ is often called \emph{completion time} in the literature, cf. \cite{takagi90}, with mean $E(B_{iL}^*) = \frac{E(B_{iL})}{1-\rho_{iH}}$. Service is exhaustive for $Q_{1H^*}, Q_{2H^*}, \dots, Q_{NH^*}$ and synchronised gated for $Q_{1L^*}, Q_{2L^*}, \dots, Q_{NL^*}$, the gate of $Q_{iL^*}$ being set at the visit beginning of $Q_{iH^*}$.
The synchronised gated service discipline is introduced in \cite{khamisy92} and does not strictly satisfy Property \ref{resingproperty}. However, it \emph{does} satisfy a slightly modified version of Property \ref{resingproperty} that still allows for straightforward analysis; see \cite{semphd} for more details. During a visit to $Q_{iL^*}$ only those type $iL^*$ customers are served that were present at the previous visit beginning to $Q_{iH^*}$. The joint queue length distribution at a visit beginning of $Q_{iH^*}$ in this system is the same as the joint queue length distribution at a visit beginning of $Q_{i}$ in the original polling system. Similarly, the joint queue length distribution at a visit completion of $Q_{iL^*}$ is the same as the joint queue length distribution at a visit completion of $Q_{i}$ in the original polling system. In terms of the GFs:
\begin{align*}
V_{b_i}(\z) &= V_{b_{iH^*}}(\z),\\
V_{c_i}(\z) &= V_{c_{iL^*}}(\z),
\end{align*}
where $\z$ is a shorthand notation for the vector $(z_{1H}, z_{1L}, \dots, z_{NH}, z_{NL})$.
The GFs of the joint queue length distributions at a visit beginning and completion of $Q_{iH^*}$ are related in the following manner:
\[
V_{c_{iH^*}}(\z) = V_{b_{iH^*}}\big(z_{1H}, z_{1L}, \dots, h_{iH}(\z), z_{iL}, \dots, z_{NH}, z_{NL}\big),
\]
with $h_{iH}(\z) = \pi_{iH}\left(\lambda_{iL}(1-z_{iL})+\sum_{j\neq i}(\lambda_{jH}(1-z_{jH}) + \lambda_{jL}(1-z_{jL}))\right)$. Similarly:
\[
V_{c_{iL^*}}(\z) = V_{b_{iH^*}}\big(z_{1H}, z_{1L}, \dots, h_{iH}(\z), h_{iL}(\z),\dots, z_{NH}, z_{NL}\big),
\]
where $h_{iL}(\z) = \beta_{iL}^*\left(\lambda_{iL}(1-z_{iL})+\sum_{j\neq i}(\lambda_{jH}(1-z_{jH}) + \lambda_{jL}(1-z_{jL}))\right)$. Note that $V_{c_{iH^*}}(\cdot) = V_{b_{iL^*}}(\cdot)$ since there is no switch-over time between $Q_{iH^*}$ and $Q_{iL^*}$. There is a switch-over time between $Q_{iL^*}$ and $Q_{(i+1)H^*}$ though:
\[
V_{b_{(i+1)H^*}}(\z) =
V_{c_{iL^*}}(\z)\sigma_i\left(\sum_{j=1}^N(\lambda_{jH}(1-z_{jH}) + \lambda_{jL}(1-z_{jL}))\right).
\]
Now that we can relate $V_{b_{(i+1)H^*}}(\cdot)$ to $V_{b_{iH^*}}(\cdot)$, we can repeat these steps $N$ times to obtain a recursive expression for $V_{b_{iH^*}}(\cdot)$. This recursive expression is sufficient to compute all moments of the joint queue length distribution at a visit beginning to $Q_{iH^*}$ by differentiation, but the expression can also be written as an infinite product which converges if and only if $\rho < 1$. We refer to \cite{resing93} for more details.

\section{Cycle time, visit time and intervisit time}\label{cycletimesection}

We define the cycle time $C_i$ as the time between two successive visit beginnings to $Q_i$, $i=1,\dots,N$. The LST of the distribution of $C_i$, denoted by $\gamma_i(\cdot)$, can be expressed in terms of $V_{b_i}(\cdot)$ because the type $iL$ customers that are present at the beginning of a visit to $Q_{i}$ are those type $iL$ customers that have arrived during the previous cycle. It is convenient to introduce the notation $\widetilde{V}_{b_i}(z_{iH}, z_{iL}) = V_{b_i}(1,\dots,1,z_{iH}, z_{iL},1,\dots,1)$, where $z_{iH}$ and $z_{iL}$ are the arguments that correspond respectively to type $iH$ and $iL$ customers. Using this notation we can write: $\widetilde{V}_{b_i}(1,z) = \gamma_i(\lambda_{iL}(1-z))$. Hence, the LST of the cycle time distribution is:
\begin{equation}
\gamma_i(\omega) = \widetilde{V}_{b_i}(1,1-\frac{\omega}{\lambda_{iL}}).   \label{cycletimelst}
\end{equation}
Note that $E(C_i) = \frac{E(S_1) +\dots+E(S_N)}{1-\rho}$, which does not depend on $i$. Higher moments of the cycle time distribution \emph{do} depend on the cycle starting point.

We define the intervisit time $I_i$ as the time between a visit completion of $Q_i$ and the next visit beginning of $Q_i$. The type $iH$ customers present at the beginning of a visit to $Q_i$ are exactly those type $iH$ customers that arrived during the previous intervisit time $I_i$. Hence, $\widetilde{V}_{b_i}(z,1) = \widetilde{I}_i(\lambda_{iH}(1-z))$, where $\widetilde{I}(\cdot)$ is the LST of the distribution of $I_i$. This leads to the following expression for the LST of the intervisit time distribution of $Q_i$:
\begin{equation}
\widetilde{I}_i(\omega) = \widetilde{V}_{b_i}(1-\frac{\omega}{\lambda_{iH}},1).  \label{intervisittimelst}
\end{equation}
It is intuitively clear that $E(I_i) = (1-\rho_{i})E(C)$.

The LSTs of the distributions of the cycle time and intervisit time are needed later in this paper. For the visit time of $Q_i$, $V_i$, we mention the LST here for completeness but it will not be used later:
\[E(\ee^{-\omega V_i}) = \widetilde{V}_{b_i}(\pi_{iH}(\omega), \beta_{iL}^*(\omega)).\]
It is easy to verify that $E(V_i) = \rho_i E(C)$.

\section{Waiting times and marginal queue lengths}\label{waitingtimelsts}

\subsection{High priority customers}

Since high priority customers are served exhaustively, we can use the concept of delay-cycles, sometimes called $T$-cycles (cf. \cite{takagi91}), introduced by Kella and Yechiali \cite{kellayechiali88} for vacation models to find the waiting time LST of a type $iH$ customer, where \emph{waiting time} is understood as the time between arrival of a customer into the system and the moment when the customer is taken into service. The waiting time plus service time will be called \emph{sojourn time} of a customer. When it comes to computing waiting times in a polling system with priorities, one can use delay-cycles for any queue that is served exhaustively, cf. \cite{boonadanboxma2queues2008,boonadanboxma2008}. A delay-cycle for a type $iH$ customer is a cycle that starts with a certain initial delay at the moment that the last type $iH$ customer in the system has been served. In our model this initial delay is either the service of a type $iL$ customer, $B_{iL}$, or (if no type $iL$ customer is present) an intervisit period $I_i$. The delay cycle ends at the first moment after the initial delay when no type $iH$ customer is present in the system again. This is the moment that all type $iH$ customers that have arrived during the delay, and all of their type $iH$ descendants, have been served. In \cite{boonadanboxma2queues2008} delay-cycles have been applied to a polling system with two priority levels in an exhaustively served queue. For a type $iH$ customer in the polling model in the present paper, the same arguments can be used to compute the LST of the waiting time distribution. The fraction of time that the system is in a delay-cycle that starts with the service time $B_{iL}$ of a type $iL$ customer is $\frac{\rho_{iL}}{1-\rho_{iH}}$, and the fraction of time that the system is in a delay-cycle that starts with an intervisit period $I_i$, is $1-\frac{\rho_{iL}}{1-\rho_{iH}} = \frac{1-\rho_i}{1-\rho_{iH}}$.
We can use the Fuhrmann-Cooper decomposition \cite{fuhrmanncooper85} to obtain the LST of the waiting time distribution of a type $iH$ customer, because from his perspective the system is an $M/G/1$ queue with server vacations. The vacation is the service time $B_{iL}$ of a type $iL$ customer with probability $\frac{\rho_{iL}}{1-\rho_{iH}}$, and an intervisit time $I_i$ with probability $\frac{1-\rho_i}{1-\rho_{iH}}$. This leads to the following expression for the LST of the waiting time distribution of a type $iH$ customer:
\begin{equation}
E[\ee^{-\omega W_{iH}}] = \frac{(1-\rho_{iH})\omega}{\omega-\lambda_{iH}(1-\beta_{iH}(\omega))}\cdot
\left[\frac{\rho_{iL}}{1-\rho_{iH}}\cdot\frac{1-\beta_{iL}(\omega)}{\omega E(B_{iL})} + \frac{1-\rho_i}{1-\rho_{iH}} \cdot \frac{1-\widetilde{I}_i(\omega)}{\omega E(I_i)}\right].\label{lstwihexhaustive}
\end{equation}
Equation \eqref{lstwihexhaustive} is similar to the equation found in \cite{boonadanboxma2queues2008} for high priority customers in an exhaustive queue. Note that the intervisit time $I_i$ is different though, with LST $\widetilde{I}_i(\cdot)$ as defined in Equation \eqref{intervisittimelst}.

The GF of the marginal queue length distribution of type $iH$ customers can be found by applying the distributional form of Little's Law \cite{keilsonservi90} to the sojourn time distribution:
\[E\left(z^{N_{iH}}\right) = E\left(\ee^{-\lambda_{iH}(1-z)(W_{iH} + B_{iH})}\right).\]
This leads to the following expression:
\begin{align}
E[z^{N_{iH}}] =& \frac{(1-\rho_{iH})(1-z)\beta_{iH}(\lambda_{iH}(1-z))}{\beta_{iH}(\lambda_{iH}(1-z))-z}\nonumber\\
&\cdot
\left[\frac{\rho_{iL}}{1-\rho_{iH}}\cdot\frac{1-\beta_{iL}(\lambda_{iH}(1-z))}{(1-z) \lambda_{iH}E(B_{iL})} + \frac{1-\rho_i}{1-\rho_{iH}} \cdot \frac{1-\widetilde{I}_i(\lambda_{iH}(1-z))}{(1-z) \lambda_{iH}E(I_i)}\right]. \label{queuelengthdecompositionH}
\end{align}

\subsection{Low priority customers}

In this subsection we determine the GF of the marginal queue length distribution of type $iL$ customers, and the LST of the waiting time distribution of type $iL$ customers. In order to obtain these functions, we regard the alternative system with $2N$ queues as defined in Section \ref{joint}. The number of type $iL$ customers in the original polling system and their waiting time (\emph{excluding} the service time) have the same distribution as the number of type $iL^*$ customers and their waiting time (again \emph{excluding} the service time, which is different) in the alternative system. From the viewpoint of a type $iL^*$ customer, the system is an ordinary polling system with synchronised gated service in $Q_{iL^*}$.

We apply the Fuhrmann-Cooper decomposition to the alternative polling model with $2N$ queues and type $iL^*$ customers having completion time $B_{iL}^*$. Using arguments similar as in the derivation of Equation (3.7) in \cite{semphd}, we find the general form of the GF of the marginal queue length distribution:
\begin{align}
E[z^{N_{iL}}] =\,& \frac{(1-\rho_{iL}^*)(1-z)\beta^*_{iL}(\lambda_{iL}(1-z))}{\beta^*_{iL}(\lambda_{iL}(1-z))-z} \nonumber\\
&\cdot
\frac{V_{c_{iL}}(1, \dots, 1, z, 1, \dots, 1) - V_{b_{iL}}(1, \dots, 1, z, 1,\dots,1)}{(1-z)(E(N^*_{iL|I_{\textit{end}}}) - E(N^*_{iL|I_{\textit{begin}}}))},   \label{queuelengthdecompositionL1}
\end{align}
where $\rho_{iL}^* = \frac{\rho_{iL}}{1-\rho_{iH}}$ and $\beta_{iL}^*(\cdot)$ is given by \eqref{extendedservicetimeL}. Furthermore, $N^*_{iL|I_{\textit{end}}}$ and $N^*_{iL|I_{\textit{begin}}}$ are the number of type $iL^*$ customers at respectively the visit \emph{beginning} and visit \emph{completion} of $Q_{iL^*}$. The visit beginning corresponds to the end of the intervisit period $I_{iL}$, and the visit completion corresponds to the beginning of the intervisit period. Substitution into \eqref{queuelengthdecompositionL1} leads to the following expression:
\begin{align}
E[z^{N_{iL}}] &= \frac{(1-\frac{\rho_{iL}}{1-\rho_{iH}})(1-z)\beta_{iL}(\lambda_{iL}(1-z)+\lambda_{iH}(1-\pi_{iH}(\lambda_{iL}(1-z))))}{\beta_{iL}(\lambda_{iL}(1-z)+\lambda_{iH}(1-\pi_{iH}(\lambda_{iL}(1-z))))-z} \nonumber\\
&\cdot
\frac{\widetilde{V}_{b_{i}}\big(\pi_{iH}(\lambda_{iL}(1-z)), \beta_{iL}(\lambda_{iL}(1-z)+\lambda_{iH}(1-\pi_{iH}(\lambda_{iL}(1-z))))\big) - \widetilde{V}_{b_{i}}\big(\pi_{iH}(\lambda_{iL}(1-z)), z\big)}{(1-z)\lambda_{iL}(1-\frac{\rho_{iL}}{1-\rho_{iH}})E(C)},\label{ezNil}
\end{align}
where we use that $E(N^*_{iL|I_{\textit{end}}}) - E(N^*_{iL|I_{\textit{begin}}}) = \lambda_{iL} (1-\rho^*_{iL})E(C) = \lambda_{iL}(1-\frac{\rho_{iL}}{1-\rho_{iH}})E(C)$, because this is the mean number of type $iL^*$ customers that arrive during the intervisit time of $Q_{iL^*}$.

Applying the distributional form of Little's Law to \eqref{ezNil}, we obtain the LST of the sojourn time distribution of a type $iL$ customer. Since the sojourn time is $W_{iL}+B^*_{iL}$, the LST of the waiting time distribution immediately follows:
\begin{align}
E[\ee^{-\omega W_{iL}}] =& \frac{(1-\frac{\rho_{iL}}{1-\rho_{iH}})\omega}{\omega-\lambda_{iL}(1-\beta_{iL}(\omega+\lambda_{iH}(1-\pi_{iH}(\omega))))} \nonumber\\
&\cdot
\frac{\widetilde{V}_{b_{i}}\big(\pi_{iH}(\omega), \beta_{iL}(\omega+\lambda_{iH}(1-\pi_{iH}(\omega)))\big) - \widetilde{V}_{b_i}\big(\pi_{iH}(\omega), 1-\frac{\omega}{\lambda_{iL}}\big)}{\omega (1-\frac{\rho_{iL}}{1-\rho_{iH}})E(C)}.\label{lstwilexhaustive}
\end{align}

\section{Moments}

Differentiation of the waiting time LSTs derived in the previous section leads to the following mean waiting times:

\begin{align}
E(W_{iH}) &= \frac{\rho_{iH} E(B_{iH,\textit{res}})+\rho_{iL} E(B_{iL,\textit{res}})}{1-\rho_{iH}}+\frac{1-\rho_i}{1-\rho_{iH}}E(I_{i,\textit{res}}),\label{EWH}\\
E(W_{iL}) &= \left(1+\frac{\rho_{iL}}{1-\rho_{iH}}\right)E(C_{i,\textit{res}}) + \frac{\rho_{iH}}{1-\rho_{iH}} \frac{E(X_{iH} X_{iL})}{\lambda_{iL} \lambda_{iH} E(C)},\label{EWL}
\end{align}
where $B_{iH,\textit{res}}$ denotes a residual service time of a type $iH$ customer, with $E(B_{iH,\textit{res}}) = \frac{E(B_{iH}^2)}{2E(B_{iH})}$. We use a similar notation for the residual service time of a type $iL$ customer, the residual intervisit time, and residual cycle time. Furthermore, $X_{iH}$ and $X_{iL}$ are respectively the number of type $iH$ and type $iL$ customers at the beginning of a visit to $Q_i$, so $E(X_{iH} X_{iL})$ is obtained by differentiating $\widetilde{V}_{b_i}(z_{iH}, z_{iL})$ with respect to $z_{iH}$ and $z_{iL}$ (and then setting $z_{iH} = z_{iL} = 1$).

We now present an alternative, direct way to obtain the mean waiting time for a type $iL$ customer by conditioning on the event that an arrival takes place in a visit period, or in an intervisit period.
\begin{align}
E(W_{iL}) =& \frac{E(V_i)}{E(C)}\left[E(V_{i,\textit{res}}) + \frac{E(V_i I_i)}{E(V_i)} + \frac{\rho_{iH}}{1-\rho_{iH}} \frac{E(V_i I_i)}{E(V_i)} + \frac{\rho_{iL}}{1-\rho_{iH}}E(V_{i,\textit{past}})\right]+\nonumber\\
&\frac{E(I_i)}{E(C)}\left[E(I_{i,\textit{res}}) + \frac{\rho_{iH}}{1-\rho_{iH}} \left(E(I_{i,\textit{past}})+E(I_{i,\textit{res}})\right) + \frac{\rho_{iL}}{1-\rho_{iH}}\left(\frac{E(V_i I_i)}{E(I_i)} + E(I_{i,\textit{past}})\right)\right]\nonumber\\
=& \frac{1}{E(C)}\left[\frac12E\left((V_i+I_i)^2\right) + \frac{\rho_{iL}}{1-\rho_{iH}}E\left((V_i+I_i)^2\right) + \frac{\rho_{iH}}{1-\rho_{iH}} \left(E(I_i^2)+E(V_i I_i)\right)\right]\nonumber\\
=& \left(1+\frac{\rho_{iL}}{1-\rho_{iH}}\right)\frac{E({C_i}^2)}{2E(C)} + \frac{\rho_{iH}}{1-\rho_{iH}} \left(\frac{E(I_i)}{E(C)} \left[2\frac{E(I_i^2)}{2E(I_i)}\right] + \frac{E(V_i)}{E(C)} \frac{E(I_i V_i)}{E(V_i)}\right).\label{EWLalternative}
\end{align}
In the above derivation, we use that both the past and residual intervisit time have expectation $\frac{E(I_i^2)}{2E(I_i)}$, and that if a type $iL$ customer arrives during the visit time (with probability $\frac{E(V_i)}{E(C)}$), the mean length of the following intervisit time equals $\frac{E(I_i V_i)}{E(V_i)}$.
The interpretation of \eqref{EWLalternative} is that a type $iL$ customer always has to wait for the residual cycle time, for the completion times of all type $iL$ customers that have arrived during the past cycle time, and for the busy periods of all type $iH$ customers that have arrived during the intervisit time of the cycle in which the type $iL$ customer has arrived.

To show that \eqref{EWL} and \eqref{EWLalternative} are equal, we can rewrite the last term in \eqref{EWL}:
\begin{align*}
E(X_{iH} X_{iL}) &= E[(N_{iL}(V_i) + N_{iL}(I_i))N_{iH}(I_i)] \\
&= E\big(E[(N_{iL}(V_i) + N_{iL}(I_i))N_{iH}(I_i)] \,|\, I_i, V_i\big) \\
&= E[(\lambda_{iL} V_i + \lambda_{iL} I_i )\lambda_{iH} I_i] \\
&= \lambda_{iL}\lambda_{iH}  E(I_i V_i)+\lambda_{iL} \lambda_{iH} E(I_i^2),
\end{align*}
where $N_{j}(T)$ denotes the number of type $j$ customers that have arrived during time $T$ $(j=iH,iL)$, and $V_i$ denotes the length of a visit of the server to $Q_i$. Hence,
\begin{align*}
\frac{E(X_{iH} X_{iL})}{\lambda_{iL} \lambda_{iH} E(C)} &= \frac{E(I_i V_i)+E(I_i^2)}{E(C)}\\
&= \frac{E(I_i)}{E(C)} \left[2\frac{E(I_i^2)}{2E(I_i)}\right] + \frac{E(V_i)}{E(C)} \frac{E(I_i V_i)}{E(V_i)},
\end{align*}
which coincides with the last term in \eqref{EWLalternative}.


\section{Pseudo-conservation law for priority polling systems}\label{pseudoconservationlawprioritiessection}

Boxma and Groenendijk \cite{boxmagroenendijk87} have shown that a so-called pseudo-conservation law holds for nonpriority polling systems. We do not discuss this law in detail in the present paper, but we mention that a generalised version of this law (cf. \cite{shimogawatakahashi88,fournierrosberg91}) holds for systems with multiple priority levels in each queue:
\begin{equation}
\begin{aligned}
\sum_{i=1}^N\sum_{k=1}^{K_i} \rho_{ik} E(W_{ik}) &=
\frac{\rho}{1-\rho} \sum_{i=1}^{N}\sum_{k=1}^{K_i} \rho_{ik}\frac{E(B_{ik}^2)}{2E(B_{ik})}\\
& + \rho \frac{E(S^2)}{2E(S)} +\left[\rho^2-\sum_{i=1}^N \rho_i^2\right]\frac{E(S)}{2(1-\rho)} + \sum_{i=1}^N E(Z_{ii}),
\end{aligned}
\label{pseudoconservationlawpriorities}
\end{equation}
where $S = \sum_{i=1}^N S_i$, and $K_i$ is the number of priority levels in $Q_i$. In this expression $Z_{ii}$ is the amount of work at $Q_i$ when the server leaves this queue and depends on the service discipline. It is well-known that for gated service, $E(Z_{ii}) = \rho_i^2 E(C)$ and for exhaustive service, $E(Z_{ii}) = 0$. The pseudo-conservation law also holds for polling systems with mixed gated/exhaustive service in some or all of the queues.
If $Q_i$ receives mixed gated/exhaustive service, we have $K_i = 2$, and $E(Z_{ii}) = \rho_{iL}\rho_i E(C)$.

\section{Numerical results}\label{numericalresults}

\subsection*{Example 1}

In order to illustrate the effect of using a mixed gated/exhaustive service discipline in a polling system with priorities, we compare it to the commonly used gated and exhaustive service disciplines. In this example we use a polling system which consists of two queues, $Q_1$ and $Q_2$. Customers in $Q_1$ are divided into high priority customers, arriving with arrival rate $\lambda_{1H} = \frac{2}{10}$, and low priority customers, with arrival rate $\lambda_{1L}=\frac{4}{10}$. Customers in $Q_2$ all have the same priority level and arrive with arrival rate $\lambda_2 = \frac{2}{10}$. All service times are exponentially distributed with mean $1$. The switch-over times $S_1$ and $S_2$ are also exponentially distributed with mean $1$, which results in a mean cycle time of $E(C) = 10$. The service discipline in $Q_2$ is gated, the service discipline in $Q_1$ is varied: gated, exhaustive and mixed gated/exhaustive. Results for a queue with two priority levels and purely gated or exhaustive service are obtained in \cite{boonadanboxma2queues2008}.

Table \ref{numericalresults1} displays the mean and the variance of the waiting times of the three customer types under the three service disciplines. We conclude that the mixed gated/exhaustive service is a major improvement for the high priority customers in $Q_1$, whereas the mean waiting times of the low priority customers in $Q_1$ and the customers in $Q_2$ hardly deteriorate. Of course in systems where $\rho_{1H}$ is quite high, the negative impact can be bigger and one has to decide exactly how far one wants to go in giving extra advantages to customers that already receive high priority. When comparing the mixed gated/exhaustive strategy to a system with purely exhaustive service in $Q_1$, we conclude that the improvement is not so much in the mean waiting time for high priority customers, but mostly in the mean and variance of the waiting time for customers in $Q_2$. 

\newcommand{\V}{\textrm{Var}}
\begin{table}[h!]
\begin{center}
\begin{tabular}{|l|r|r|r|}
\hline
& Gated & Exhaustive & Mixed G/E \\
\hline
$E(W_{1H})$        &    9.578 & 2.520 & 2.338 \\
$E(W_{1L})$        &    14.366 & 6.300 & 14.575 \\
$E(W_{2})$         &    9.690 & 14.880 & 10.513 \\
\hline$\V(W_{1H})$ &    56.739 & 9.290 & 6.496 \\
$\V(W_{1L})$       &    101.616 & 32.812 & 118.217 \\
$\V(W_{2})$        &    58.513 & 231.256 & 76.371
\\\hline
\end{tabular}
\caption{Numerical results for Example 1. The switch-over times $S_1$ and $S_2$ are exponentially distributed with mean 1. The mixed gated/exhaustive service discipline is compared to gated and exhaustive service.}
\label{numericalresults1}
\end{center}
\end{table}

\begin{table}[h!]
\begin{center}
\begin{tabular}{|l|r|r|r|}
\hline
& Gated & Exhaustive & Mixed G/E \\
\hline
$E(W_{1H})$        &     63.187 & 11.333 & 11.167 \\
$E(W_{1L})$        &     94.781 & 28.333 & 90.417 \\
$E(W_{2})$         &     63.251 & 68.000 & 64.000 \\
\hline$\V(W_{1H})$ &     847.377 & 195.508 & 183.907 \\
$\V(W_{1L})$       &     894.173 & 315.823 & 850.199 \\
$\V(W_{2})$        &     853.777 & 1386.100 & 928.914
\\\hline
\end{tabular}
\caption{Numerical results for Example 1. Switch-over times are deterministic: $S_1 = S_2 = 10$. 
}
\label{numericalresults2}
\end{center}
\end{table}

It is noteworthy that the mixed gated/exhaustive service discipline does not always have a negative effect on the mean waiting time of low priority customers in $Q_1$, $E(W_{1L})$, compared to the gated service discipline. If, for example, the switch-over times are taken to be deterministic with value 10, the mean waiting time for low priority customers is significantly less for the mixed gated/exhaustive service than for gated service, as can be seen in Table \ref{numericalresults2}. Compared to gated service, type $1H$ customers benefit strongly from the mixed gated/exhaustive service discipline, and even type $1L$ customers benefit from it. The mean waiting time for customers in $Q_2$ has increased, but only marginally.

In order to get more understanding of this surprising behaviour of the waiting time of low priority customers as function of the arrival intensities $\lambda_{1H}$ and $\lambda_{1L}$, we use a simplified model which leads to more insightful expressions, but displays the same characteristics as the model that was analysed in the previous paragraph. Instead of analysing a polling model, we analyse an $M/G/1$ queue with multiple server vacations. The queue, denoted by $Q_1$ to use familiar notation, contains high (type $1H$) and low (type $1L$) priority customers. Also here high priority customers are served before low priority customers. The service times of both customers types are exponentially distributed with mean 1. This is for notational reasons only, for this example we actually only require that both service times are identically distributed. One server vacation has a fixed length $S$. If the server does not find any customers waiting upon arrival from a vacation, he takes another vacation of length $S$ and so on. In order to stay consistent with the notation used earlier, we denote the occupation rate of high and low priority customers by respectively $\rho_{1H}$ and $\rho_{1L}$. The total occupation rate is $\rho = \rho_1 = \rho_{1H} + \rho_{1L}$. Note that in this example $\lambda_{1H} = \rho_{1H}$ and $\lambda_{1L} = \rho_{1L}$. We now compare the mean waiting times of type $1L$ customers in the system with purely gated service and the system with mixed gated/exhaustive service. For this simplified model, we can write down explicit expressions that have been obtained by differentiating the LSTs and solving the resulting equations. These expressions could also have been obtained by using Mean Value Analysis (MVA) for polling systems \cite{wierman07,winands06}.

\begin{align}
&\textrm{Gated service:} & &  E(W_{1L}) = (1 + \rho + \rho_{1H})\left(\frac{S}{2(1-\rho)} + \frac{\rho}{1-\rho^2}\right),\label{EWLgatedvacation}\\
&\textrm{Mixed G/E service: } & &  E(W_{1L}) = \frac{\rho}{(1-\rho)(1-\rho_{1H})} + \frac{S(1+\rho(1-2\rho_{1H}))}{2(1-\rho)(1-\rho_{1H})}.\label{EWLgatedexhaustivevacation}
\end{align}

Now we analyse the behaviour of these waiting times as we vary $\lambda_{1H}$ between 0 and $\rho$, while keeping $\lambda_{1H}+\lambda_{1L} = \rho$ constant. Substitution of $\lambda_{1H} = 0$ shows that the mean waiting times in the gated and mixed gated/exhaustive system are equal:
\[E(W_{1L}|\rho_{1H}=0) = \frac{S(1 + \rho)}{2(1-\rho)} + \frac{\rho}{1-\rho}.\]
Letting $\lambda_{1H} \rightarrow \rho$ leads to the following expressions:

\begin{tabular}{ll}
Gated service:  & $\displaystyle E(W_{1L}|\rho_{1H} \rightarrow \rho) = \frac{\rho(1 + 2\rho)}{1-\rho^2}+\frac{S(1 + 2\rho)}{2(1-\rho)}$,\\[3ex]
Mixed G/E service:  & $\displaystyle E(W_{1L}|\rho_{1H} \rightarrow \rho) = \frac{\rho}{(1-\rho)^2} + \frac{S(1+2\rho)}{2(1-\rho)}$.
\end{tabular}

Two interesting things can be concluded from these two equations for the case $\lambda_{1H} \rightarrow \rho$:
\begin{itemize}
\item for fixed $\rho$, $E(W_{1L})$ in a gated system is always less than $E(W_{1L})$ in a mixed gated/exhaustive system,
\item the difference between $E(W_{1L})$ in a gated system and $E(W_{1L})$ in a mixed gated/exhaustive system does \emph{not} depend on $S$.
\end{itemize}

Focussing on the mean waiting time of type $1L$ customers only, we conclude that a gated system performs the same as a mixed gated/exhaustive system as $\rho_{1L} = \rho$, and that a gated system always performs better when $\rho_{1L} \rightarrow 0$. For $0 < \rho_{1L} < \rho$ the vacation time $S$ determines which system performs better. By taking derivatives of \eqref{EWLgatedvacation} and \eqref{EWLgatedexhaustivevacation} with respect to $\rho_{1H}$ and letting $\rho_{1H} \rightarrow 0$, one finds that the mean waiting time of a type $1L$ customer in a mixed gated/exhaustive system is less than in a purely gated system when $\rho_{1H} \rightarrow 0$, if and only if $S > \frac{2\rho}{1+\rho}$. Since a gated system always outperforms a mixed gated/exhaustive system when $\lambda_{1H} \rightarrow \rho$, for $S > \frac{2\rho}{1+\rho}$ there must be (at least) one value of $\lambda_{1H}$ for which the two systems perform the same. Further inspection of the derivatives gives the insight that in a gated system the relation between $E(W_{1L})$ and $\lambda_{1H}$ is a straight line, which can also be seen immediately from Equation \eqref{EWLgatedvacation}. In a mixed gated/exhaustive system, the relation between $E(W_{1L})$ and $\lambda_{1H}$ is not a straight line, both the first and second derivative with respect to $\lambda_{1H}$ are strictly positive. This means that for $S \leq \frac{2\rho}{1+\rho}$ the gated system always performs better than the mixed gated/exhaustive system for any value of $\lambda_{1H} > 0$, and for $S > \frac{2\rho}{1+\rho}$ the mixed gated/exhaustive system performs better than the gated system for $0 < \lambda_{1H} < \lambda_{1H}^*$. The value of $\lambda_{1H}^*$ can be determined analytically:
\[\lambda_{1H}^* = \rho\frac{S-\frac{2\rho}{1+\rho}}{S+\frac{2\rho}{1+\rho}}.\]
From this expression we conclude that $\lim_{S \rightarrow \infty} \lambda_{1H}^* = \rho$. Although we have studied only the vacation model, the conclusions are also valid for more general settings, like polling models with non-deterministic switch-over times, but the expressions are by far not as appealing. 

\begin{figure}[h!]
\begin{center}
\includegraphics[width=0.2\textwidth]{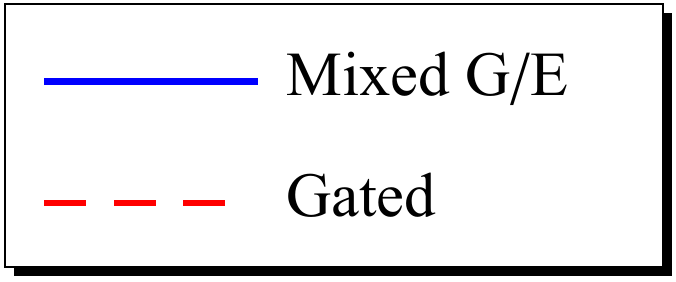}\\
\parbox{0.32\textwidth}{\centering
\includegraphics[width=\linewidth]{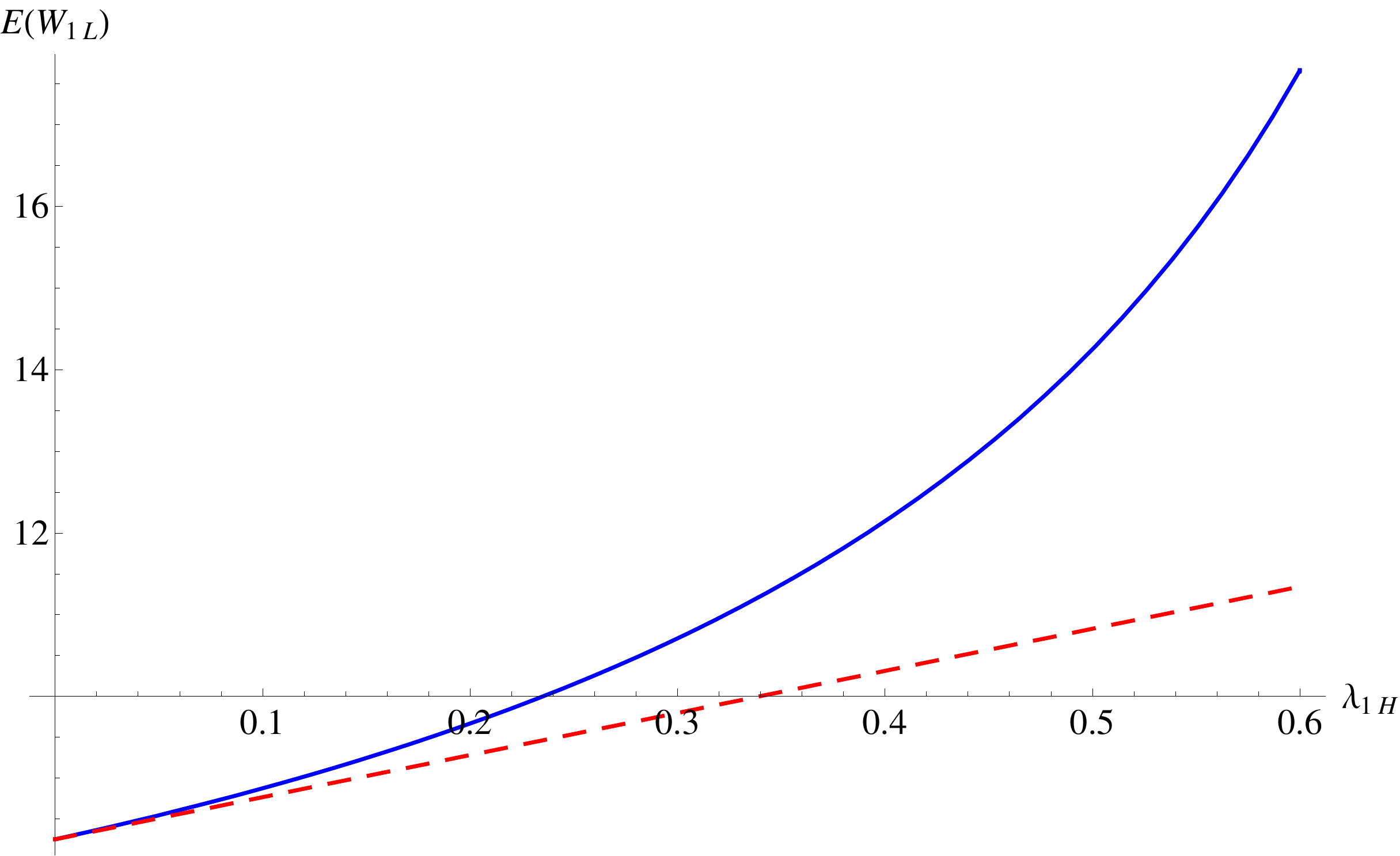}\\a) $S=1$
}
\hfill
\parbox{0.32\textwidth}{\centering
\includegraphics[width=\linewidth]{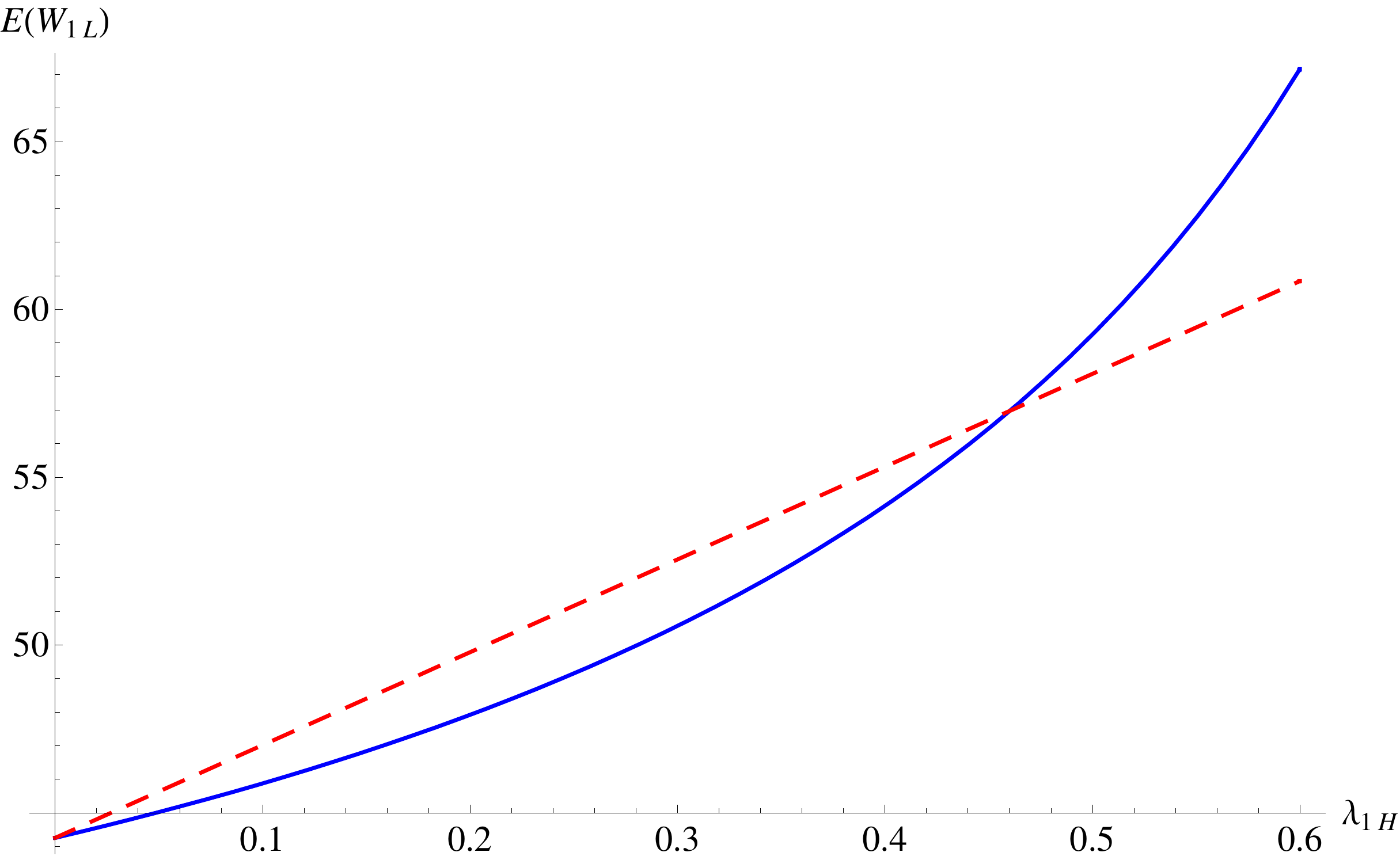}\\b) $S=10$
}
\hfill
\parbox{0.32\textwidth}{\centering
\includegraphics[width=\linewidth]{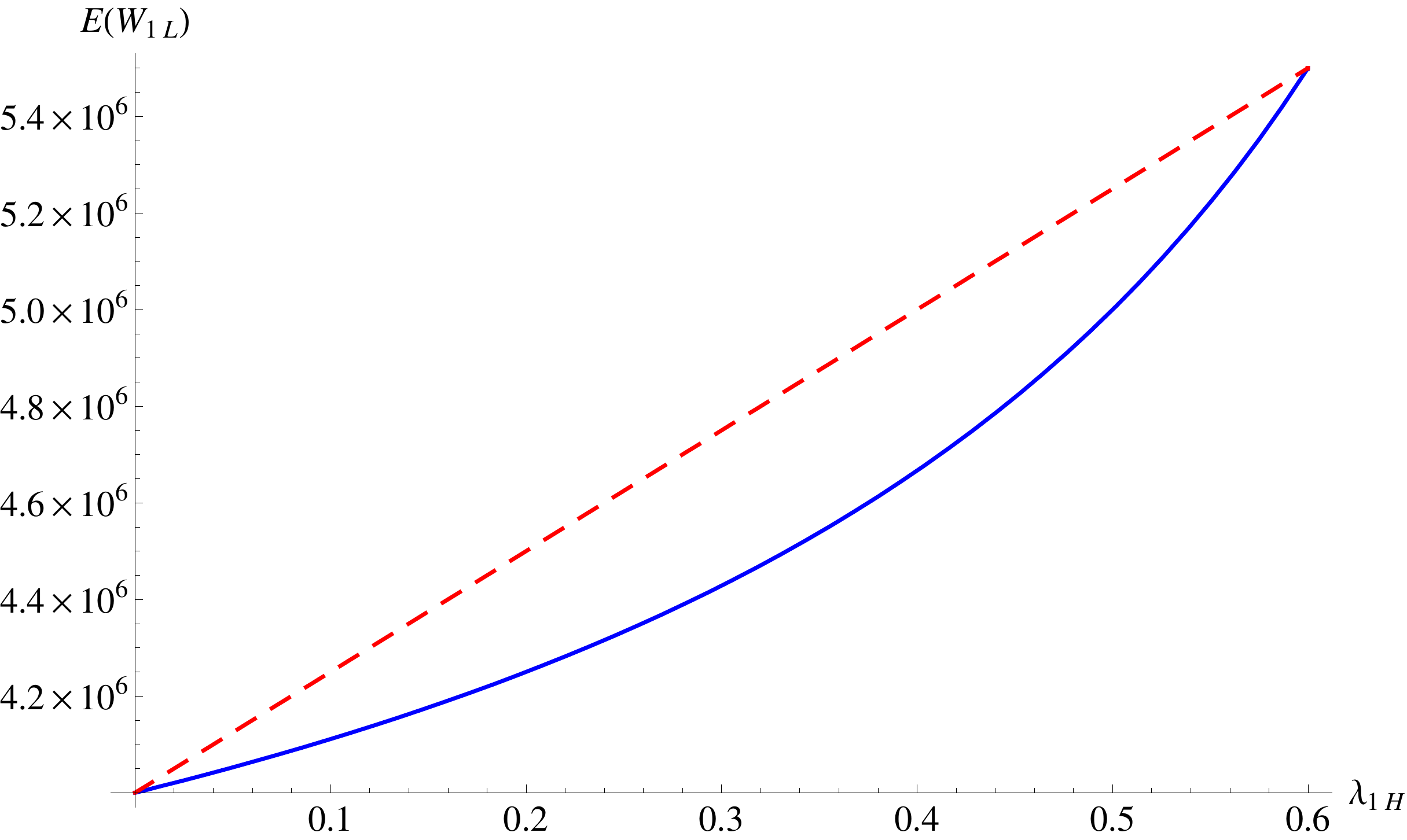}\\c) $S=10^6$
}
\end{center}
\caption{Mean waiting time of type $1L$ customers in the polling model discussed in Example 1. For gated and mixed gated/exhaustive service $E(W_{1L})$ is plotted against $\lambda_{1H}$ while keeping $\lambda_{1L} + \lambda_{1H}$ constant. The switch-over times $S_1 = S_2 = S/2$ are deterministic.}
\label{ew1lGatedVsMixedGE}
\end{figure}

We visualise the findings of the present section in Figure \ref{ew1lGatedVsMixedGE}, where we show three plots of the mean waiting time of type $1L$ customers against $\lambda_{1H}$. The model considered is the same as in the beginning of the present section (two queues, gated service in $Q_2$) except for the switch-over times $S_1$ and $S_2$, which are now deterministic. We compare gated service in $Q_1$ to mixed gated/exhaustive service for three different switch-over times (notice that the scales of the three plots in Figure \ref{ew1lGatedVsMixedGE} are different).

\subsection*{Example 2}

In the previous example we showed that the mixed gated/exhaustive service discipline does not necessarily have a negative impact on the mean waiting times of low priority customers. In this example we aim at giving a better comparison of the performance of the gated, exhaustive and mixed gated/exhaustive service disciplines in a polling system with priorities. The polling system considered consists of two queues, each having high and low priority customers. The switch-over times $S_1$ and $S_2$ are exponentially distributed with mean 10. Service times of all customer types are exponentially distributed with mean 1. The arrival rates of the various customer types are: $\lambda_{1H} = \lambda_{1L} = \frac{1}{10}$, and $\lambda_{2H} = \lambda_{2L} = \frac{7}{20}$. The total occupation rate of this polling system is $\rho = \frac{9}{10}$, and we deliberately choose a system where the occupation rates of the two queues are very different, and the switch-over times are relatively high compared to the service times. The reason is that we envision production systems as the main application for the present paper (see also Section \ref{applicationssection}). In these applications large setup times are very common (see, e.g., \cite{winandsPhD}).

Table \ref{tableexample2} shows the mean and variance of the waiting times of all customer types of this polling system for all combinations of gated, exhaustive and mixed gated/exhaustive service. 
We leave it up to the reader to pick his favourite combination of service disciplines, but our preference goes out to the system with exhaustive service in $Q_1$ and mixed gated/exhaustive service in $Q_2$ because in our opinion the best combination of low mean waiting times and moderate variances is obtained in this system.

\begin{table}[h!]
\setlength{\parskip}{3ex}
\begin{center}
\begin{tabular}{|l|l|rr|rr|}
\hline
Queue& Service discipline& $E(W_{iL})$ & $E(W_{iH})$ & $\V(W_{iL})$ & $\V(W_{iH})$ \\
\hline
 1 & \text{gated} & 141.81 & 119.99 & 5166.03 & 4660.09 \\
 2 & \text{gated} & 222.95 & 146.82 & 5917.70 & 3560.67
\\\hline
\end{tabular}

\begin{tabular}{|l|l|rr|rr|}
\hline
Queue& Service discipline& $E(W_{iL})$ & $E(W_{iH})$ & $\V(W_{iL})$ & $\V(W_{iH})$ \\
\hline
 1 & \text{gated} & 165.49 & 140.03 & 11087.40 & 9411.43 \\
 2 & \text{exhaustive} & 59.45 & 17.83 & 1862.57 & 651.03
\\\hline
\end{tabular}

\begin{tabular}{|l|l|rr|rr|}
\hline
Queue& Service discipline& $E(W_{iL})$ & $E(W_{iH})$ & $\V(W_{iL})$ & $\V(W_{iH})$ \\
\hline
 1 & \text{gated} & 147.38 & 124.71 & 6406.11 & 5658.44 \\
 2 & \text{gatedexhaustive} & 209.86 & 16.98 & 6213.92 & 555.67
\\\hline
\end{tabular}

\begin{tabular}{|l|l|rr|rr|}
\hline
Queue& Service discipline& $E(W_{iL})$ & $E(W_{iH})$ & $\V(W_{iL})$ & $\V(W_{iH})$ \\
\hline
 1 & \text{exhaustive} & 97.63 & 78.10 & 4252.19 & 3784.99 \\
 2 & \text{gated} & 224.00 & 147.51 & 6186.88 & 3690.81
\\\hline
\end{tabular}

\begin{tabular}{|l|l|rr|rr|}
\hline
Queue& Service discipline& $E(W_{iL})$ & $E(W_{iH})$ & $\V(W_{iL})$ & $\V(W_{iH})$ \\
\hline
 1 & \text{exhaustive} & 119.80 & 95.84 & 9516.58 & 7952.09 \\
 2 & \text{exhaustive} & 61.62 & 18.49 & 2136.19 & 728.97
\\\hline
\end{tabular}

\begin{tabular}{|l|l|rr|rr|}
\hline
Queue& Service discipline& $E(W_{iL})$ & $E(W_{iH})$ & $\V(W_{iL})$ & $\V(W_{iH})$ \\
\hline
 1 & \text{exhaustive} & 102.18 & 81.75 & 5193.21 & 4533.58 \\
 2 & \text{gatedexhaustive} & 211.90 & 17.27 & 6722.53 & 586.84
\\\hline
\end{tabular}

\begin{tabular}{|l|l|rr|rr|}
\hline
Queue& Service discipline& $E(W_{iL})$ & $E(W_{iH})$ & $\V(W_{iL})$ & $\V(W_{iH})$ \\
\hline
 1 & \text{gatedexhaustive} & 140.95 & 77.96 & 5140.20 & 3756.12 \\
 2 & \text{gated} & 223.45 & 147.15 & 6045.55 & 3622.49
\\\hline
\end{tabular}

\begin{tabular}{|l|l|rr|rr|}
\hline
Queue& Service discipline& $E(W_{iL})$ & $E(W_{iH})$ & $\V(W_{iL})$ & $\V(W_{iH})$ \\
\hline
 1 & \text{gatedexhaustive} & 166.85 & 94.38 & 11655.90 & 7574.67 \\
 2 & \text{exhaustive} & 60.39 & 18.12 & 1978.87 & 684.25
\\\hline
\end{tabular}

\begin{tabular}{|l|l|rr|rr|}
\hline
Queue& Service discipline& $E(W_{iL})$ & $E(W_{iH})$ & $\V(W_{iL})$ & $\V(W_{iH})$ \\
\hline
 1 & \text{gatedexhaustive} & 146.87 & 81.41 & 6452.48 & 4462.04 \\
 2 & \text{gatedexhaustive} & 210.82 & 17.10 & 6451.10 & 569.08
\\\hline
\end{tabular}
\end{center}
\caption{Expectation and variance of the waiting times of the polling model discussed in Section \ref{numericalresults}, Example 2.}
\label{tableexample2}
\end{table}

\section{Possible extensions and variations}\label{extensions}

Many extensions or variations of the model discussed in the present paper can be thought of. In this section we discuss some of them.

\paragraph{A globally gated system.} The globally gated service discipline has received quite some attention in polling systems. Instead of setting the gates at the beginning of a visit to a certain queue, the globally gated service discipline states that all gates are set at the beginning of a cycle, which is the start of a visit to an arbitrarily chosen queue. The model under consideration can be analysed using similar techniques if high priority customers are served exhaustively, but low priority customers are served according to the globally gated service discipline. One would first have to build a similar model that contains $2N$ queues and determine the joint queue length distribution at visit beginnings and endings. The cycle time, starting at the moment that all gates are set, can be expressed in terms of the GF of the number of customers at the beginning of that cycle. Waiting times for high priority customers can be obtained using delay-cycles again, and waiting times for low priority customers can be obtained using the Fuhrmann-Cooper decomposition.
The LST of the waiting time distribution of low priority customers gets more complicated as the queue gets served later in the cycle.

\paragraph{More than two priority levels.} It is possible to analyse a similar model as the one of Section \ref{modeldescription}, but with more than two, say $K_i$, priority levels in $Q_i$. These $K_i$ priority levels still have to be divided into two categories: high priority levels $1, \dots, k_i$ that receive exhaustive service, and low priority levels $k_i+1, \dots, K_i$ that receive gated service. The methodology from Section \ref{waitingtimelsts} can be used, combined with the techniques that are used to analyse a polling model with multiple priority levels, cf. \cite{boonadanboxma2008}.

\paragraph{A mixture of gated and exhaustive without priorities.} One could think of a system where each queue contains two customer classes having respectively the exhaustive and gated service discipline, but service is First-Come-First-Served (FCFS). The model is similar to the model discussed in this paper, with the exception that no ``overtaking'' takes place. Customers that are served exhaustively will not be served before any ``gated customers'' standing in front of this gate, but they are allowed to pass the gate. The joint queue length distributions at polling epochs and the cycle times are the same as for the system considered in the present paper. Since no overtaking takes place, the waiting times can be found without the use of delay cycles. Nevertheless, analysis of the waiting times is quite tedious because a visit of a server to $Q_i$ consists of three parts. The third part is the service of exhaustive customers behind the gate, the first part is the service of the gated customers that have arrived during the ``previous third part'' and the second part is the FCFS service of both gated and exhaustive customers that have arrived during the previous intervisit time of $Q_i$. A combination of this non-priority mixture of gated and exhaustive, and the service discipline discussed in the present paper is discussed by Fiems et al. \cite{fiemsdevuystbruneel02}. They introduce, albeit in the different setting of a vacation queue modelled in discrete time, a service discipline where high priority customers in front of the gate are served before low priority customers waiting in front of the gate. The difference with the model discussed in the present paper, is that high priority customers entering the queue while it is being visited can pass the gate, but are not allowed to overtake low priority customers standing in front of the gate. 

\section*{Acknowledgements}

The authors wish to thank Erik Winands for his many helpful remarks and discussions. His contribution to the present paper is very much appreciated. Our gratitude also goes out to Jacques Resing who suggested the mixed gated/exhaustive service discipline. Finally, the authors thank Onno Boxma for valuable discussions and for useful comments on earlier drafts of the present paper.

\bibliographystyle{abbrvnat}

\end{document}